# The Natural Lift Curves for the Spherical Indicatrices of the Timelike-Spacelike Bertrand Couple in Minkowski 3-Space

M. BİLİCİ[1*] E. ERGÜN[2] and M. ÇALIŞKAN[3]

[1]Ondokuz Mayıs University, Faculty of Arts and Sciences, Department of Mathematics,
55200 Kurupelit, Samsun, TURKEY, E-mail: mbilici@omu.edu.tr

[2]Ondokuz Mayıs University, Çarşamba Chamber of Commerce Vocational School, 55500 Çarşamba,
Samsun,Turkey, E-mail: eergun@omu.edu.tr

[3*]Gazi University, Faculty of Sciences, Department of Mathematics, 06500 Teknik Okullar,
Ankara, Turkey, E-mail: mustafacaliskan@gazi.edu.tr

**Abstract** : This paper deals with the natural lift curves and the geodesic sprays for the spherical indicatrices of the timelike-spacelike Bertrand couple on the tangent bundle $T(S_1^2)$ or $T(H_0^2)$ in Minkowski 3-space and then give some new characterizations for these curves. Additionally we illustrate an example of our main results.

## 1. Introduction

The notion of Bertrand curves was discovered by J. Bertrand in 1850 and it plays an important role in classical differential geometry. A Bertrand curve is a curve such that its principal normals are the principal normals of a second curve (called Bertrand mate). It is proved in most texts on the subject that the characteristic property of such a curve is the existence of a linear relation between the curvature and the torsion; the discussion appears as an application of the Frenet-Serret formulas. So, a circular helix is a Bertrand curve. Bertrand mates represent particular examples of offset curves [14] which are used in computer-aided design (CAD) and computer-aided manufacturing (CAM). For classical and basic treatments

---
[*] Corresponding author

of Bertrand curves, we refer to [6], [10], [18] and [12].

There are recent works about the Bertrand curves. Ekmekçi and İlarslan studied Nonnull Bertrand curves in the n-dimensional Lorentzian space. Straightforward modication of classical theory to spacelike or timelike curves in Minkowski 3-space is easily obtained, (see [7]). Izumiya and Takeuchi [8] have shown that cylindrical helices can be constructed from plane curves and Bertrand curves can be constructed from spherical curves. Also, the representation formulae for Bertrand curves were given by [16].

In differential geometry, especially the theory of space curves, the Darboux vector is the areal velocity vector of the Frenet frame of a space curve. It is named after Gaston Darboux who discovered it. In terms of the Frenet-Serret apparatus, the Darboux vector $w$ can be expressed as $w = \tau t + \kappa b$. In addition, the concepts of the natural lift and the geodesic sprays have first been given by Thorpe (1979). On the other hand, Çalışkan et al. [5] have studied the natural lift curves and the geodesic sprays in Euclidean 3-space $\mathbb{R}^3$. Bilici et al. [2] have proposed the natural lift curves and the geodesic sprays for the spherical indicatrices of the involute-evolute curve couple in $\mathbb{R}^3$. Then, Bilici [3] adapted this problem for the spherical indicatrices of the involutes of a timelike curve in Minkowski 3-space. However, this problem is not handled in other types of space curves.

Spherical images (indicatrices) are a well known concept in classical differential geometry of curves [13]. Kula and Yaylı [11] have studied spherical images of the tangent indicatrix and binormal indicatrix of a slant helix and they have shown that the spherical images are spherical helices. In [21] Yılmaz et. al. investigated tangent and trinormal spherical images of timelike curve lying on the pseudo hyperbolic space $H_0^3$ in Minkowski space-time. İyigün [9] defined the tangent spherical image of a unit speed timelike curve lying on the on the pseudo hyperbolic space $H_0^2$ in $\mathbb{R}_1^3$.

In this study, we carry tangents of the spacelike Bertrand mate $\beta$ of a timelike curve $\alpha$ to the center of the unit hypersphere $S_1^2$ and we obtain a spacelike curve $\beta_{T^*} = T^*$ on the unit hypersphere $S_1^2$. This curve is called the first spherical indicatrix or tangent indicatrix of $\beta$. One consider the principal normal indicatrix $\beta_{N^*} = N^*$ and the binormal indicatrix $\beta_{B^*} = B^*$ on the unit hypersphere $S_1^2$ and $H_0^2$, respectively. Then the natural lift curves and the geodesic sprays for the spherical indicatrices of the timelike-spacelike Bertrand couple are investigated in Minkowski 3-space and some new results are obtained.

## 2. Preliminaries

To meet the requirements in the next sections, the basic elements of the theory of curves and hypersurfaces in the Minkowski 3-space are briefly presented in this section. A more detailed information can be found in [15].

The Minkowski 3-space $\mathbb{R}_1^3$ is the real vector space $\mathbb{R}^3$ endowed with standard flat Lorentzian metric given by

$$g = -dx_1^2 + dx_2^2 + dx_3^2,$$

where $(x_1, x_2, x_3)$ is a rectangular coordinate system of $\mathbb{R}_1^3$. A vector $V = (v_1, v_2, v_3) \in IR_1^3$ is said to be timelike if $g(V,V) < 0$, spacelike if $g(V,V) > 0$ or $V = 0$ and null (lightlike) if $g(V,V) = 0$ or $V \neq 0$. Similarly, an arbitrary curve $\Gamma = \Gamma(s)$ in $\mathbb{R}_1^3$ can locally be timelike, spacelike or null (lightlike), if all of its velocity vectors $\Gamma'$ are respectively timelike, spacelike or null (lightlike), for every $t \in I \subset \mathbb{R}$. The pseudo-norm of an arbitrary vector $V \in \mathbb{R}_1^3$ is given by $\|V\| = \sqrt{|g(V,V)|}$. $\Gamma$ is called a unit speed curve if the velocity vector $V$ of $\Gamma$ satisfies $\|V\| = 1$. A timelike vector $V$ is said to be positive (resp. negative) if and only if $v_1 > 0$ (resp. $v_1 < 0$). Let $\Gamma$ be a unit speed timelike curve with curvature $\kappa$ and torsion $\tau$. Denote by $\{t(s), n(s), b(s)\}$ the moving Frenet frame along the curve $\Gamma$ in the space $\mathbb{R}_1^3$. Then $t$, $n$ and $b$ are the tangent, the principal normal and the binormal vector of the curve $\Gamma$, respectively. For these vectors, we can write

$$t \times n = -b, \quad n \times b = t, \quad b \times t = -n,$$

where '$\times$' is the Lorentzian cross product in space $\mathbb{R}_1^3$. Depending on the causal character of the curve $\Gamma$, the following Frenet formulae are given in [17].

$$\begin{cases} \dot{t} = \kappa n, \quad \dot{n} = \kappa t - \tau b, \quad \dot{b} = \tau n \\ g(t,t) = -1, \; g(n,n) = g(b,b) = 1, \; g(t,n) = g(t,b) = g(n,b) = 0 \end{cases}$$

The Darboux vector for the timelike curve is defined by [20]:

$$w = \tau t - \kappa b.$$

There are two cases corresponding to the causal characteristic of Darboux vector $w$.

**Case I.** If $|\kappa| > |\tau|$, then $w$ is a spacelike vector. In this situation, we can write

$$\begin{cases} \kappa = \|w\| \cosh\varphi \\ \tau = \|w\| \sinh\varphi \end{cases}, \quad \|w\|^2 = g(w,w) = \kappa^2 - \tau^2$$

and the unit vector $c$ of direction $w$ is

$$c = \frac{1}{\|w\|} w = \sinh\varphi \, t - \cosh\varphi \, b,$$

where $\theta$ is the lorentzian timelike angle between $-b$ and timelike unit vector $c'$, Lorentz orthogonal to the normalisation of the Darboux vector $c = \frac{1}{\|w\|} w$ as Figure 1.

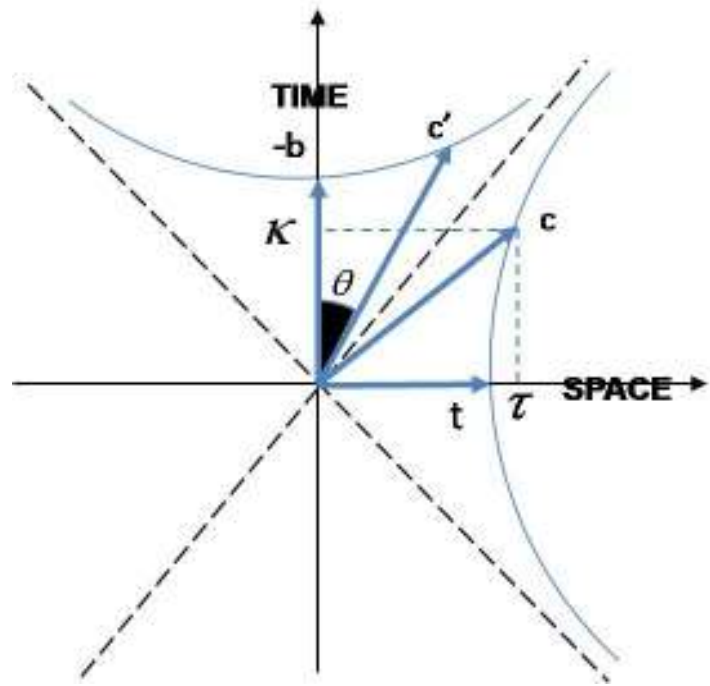

**Figure 1:** Lorentzian timelike angle $\theta$

**Case II.** If $|\kappa| < |\tau|$, then $w$ is a timelike vector. In this situation, we have

$$\begin{cases} \kappa = \|w\| \sinh\varphi \\ \tau = \|w\| \cosh\varphi \end{cases}, \quad \|w\|^2 = -g(w,w) = \tau^2 - \kappa^2$$

and the unit vector $c$ of direction $w$ is

$$c = \frac{1}{\|w\|} w = \sinh\varphi \, t - \cosh\varphi \, b.$$

**Proposition 2.1.** Let $\alpha$ be a timelike (or spacelike) curve with curvatures $\kappa$ and $\tau$. The curve $\alpha$ is a general helix if and only if $\dfrac{\tau}{\kappa} = costant$, [1].

**Remark 2.1.** We can easily see from equations of the section *Case I* and *Case II* that: $\dfrac{\tau}{\kappa} = tanh\varphi$ or $\dfrac{\tau}{\kappa} = coth\varphi$, if $\varphi = constant$ then $\alpha$ is a general helix.

**Lemma 2.1.** *The natural lift $\bar{\alpha}$ of the curve $\alpha$ is an integral curve of the geodesic spray X if and only if $\alpha$ is a geodesic on M*, [3].

**Definition 2.1.** Let $\alpha = (\alpha(s); T(s), N(s), B(s))$ be a regular timelike curve and $\beta = (\beta(s^*); T^*(s^*), N^*(s^*), B^*(s^*))$ be a regular spacelike curve in $\mathbb{R}^3_1$. $\alpha(s)$ and $\beta(s^*)$ are called Bertrand curves if $N(s)$ and $N^*(s^*)$ are linearly dependent. In this situation, $(\alpha, \beta)$ is called a timelike-spacelike Bertrand couple. (See [7] for the more details in the n-dimensional space)

**Lemma 2.2.** Let $(\alpha, \beta)$ be a timelike-spacelike Bertrand couple. The relations between the Frenet vectors of the $(\alpha, \beta)$ is as follow

$$\begin{bmatrix} T^* \\ N^* \\ B^* \end{bmatrix} = \begin{bmatrix} \sinh\theta & 0 & \cosh\theta \\ 0 & 1 & 0 \\ \cosh\theta & 0 & \sinh\theta \end{bmatrix} \begin{bmatrix} T \\ N \\ B \end{bmatrix},$$

where $\theta$ is the hyperbolic timelike angle between $T$ and $T^*$, [16].

**Remark 2.2.** The measure of the angle between the vector fields of Bertrand curve is constant, $g(T, T^*) = -\sinh\theta = constant$.

**Definition 2.2.** Let $S^2_1$ and $H^2_0$ be hypersphere in $\mathbb{R}^3_1$. The Lorentzian sphere and hyperbolic sphere of radius 1 in $IR^3_1$ are given by

$$S^2_1 = \{V = (v_1, v_2, v_3) \in \mathbb{R}^3_1 : g(V,V) = 1\}$$

and

$$H^2_0 = \{V = (v_1, v_2, v_3) \in \mathbb{R}^3_1 : g(V,V) = -1\}$$

respectively, [15].

**Definition 2.3.** Let $M$ be a hypersurface in $\mathbb{R}^3_1$ equipped with a metric $g$. Let $TM$ be the set $\cup\{T_p(M): p \in M\}$ of all tangent vectors to $M$. Then each $v \in TM$ is in a unique $T_p(M)$, and the projection $\pi : TM \to M$ sends $v$ to $p$. Thus $\pi^{-1}(p) = T_p(M)$. There is a natural way to make $TM$ a manifold, called the *tangent bundle* of $M$.

A vector field $X \in \chi(M)$ is exactly a smooth section of $TM$, that is, a smooth function $X : M \to TM$ such that $\pi \circ X = id_M$, [3].

**Definition 2.4.** Let $M$ be a hypersurface in $\mathbb{R}_1^3$. A curve $\alpha : I \to M$ is an *integral curve* of $X \in \chi(M)$ provided $\dot{\alpha} = X_\alpha$; that is,

$$\frac{d}{dt}(\alpha(s)) = X(\alpha(s)) \text{ for all } s \in I, [15]. \tag{1}$$

**Definition 2.5.** For any parametrized curve $\alpha : I \to M$, the parametrized curve

$$\overline{\alpha} : I \to TM$$

given by

$$\overline{\alpha}(t) = (\alpha(s), \dot{\alpha}(s)) = \dot{\alpha}(s)|_{\alpha(s)} \tag{2}$$

is called the *natural lift* of $\alpha$ on $TM$. Thus, we can write

$$\frac{d\overline{\alpha}}{dt} = \frac{d}{dt}(\dot{\alpha}(s))|_{\alpha(s)} = D_{\dot{\alpha}(s)}\dot{\alpha}(s) \tag{3}$$

where $D$ is the standart connection on $\mathbb{R}_1^3$, [3].

**Definition 2.6.** For $v \in TM$, a smooth vector field $X \in \chi(TM)$ defined by

$$X(v) = \varepsilon\, g(v, S(v))\xi|_{\alpha(t)}, \quad \varepsilon = g(\xi, \xi) \tag{4}$$

is called the *geodesic spray* on the manifold $TM$, where $\xi$ is the unit normal vector field of $M$ and $S$ is the shape operator of $M$, [3].

## 3. The Natural Lift Curves of the Spherical Indicatrices of the Timelike-Spacelike Bertrand Couple in Minkowski 3-Space

In this section we investigate the natural lift curves of the spherical indicatrices of spacelike Bertrand mate $\beta$ of a timelike curve $\alpha$. Furthermore, some interesting theorems about the original curve were obtained depending on the assumption that the natural lift curves should be the integral curve of the geodesic spray on the tangent bundle $T(S_1^2)$ or $T(H_0^2)$.

Note that $\overline{D}$ and $\overline{\overline{D}}$ are Levi-Civita connections on $S_1^2$ and $H_0^2$, respectively. Then Gauss equations are given by the followings

$$D_X Y = \overline{D}_X Y + \varepsilon g(S(X), Y)\xi, \quad D_X Y = \overline{\overline{D}}_X Y + \varepsilon g(S(X), Y)\xi, \quad \varepsilon = g(\xi, \xi),$$

where $\xi$ is a unit normal vector field and $S$ is the shape operator of $S_1^2$ (or $H_0^2$).

### 3.1. The natural lift of the spherical indicatrix of the tangent vector of β

Let $\alpha$ be a timelike curve and $\beta$ be a spacelike Bertrand mate for $\alpha$. We will investigate the curve $\alpha$ to satisfy the condition that the natural lift curve of $\overline{\beta}_{T^*}$ is an integral curve of geodesic spray, where $\beta_{T^*}$ is the tangent indicatrix of $\beta$. If the natural lift curve $\overline{\beta}_{T^*}$ is an integral curve of the geodesic spray, then by means of Lemma 2.1. we get,

$$\overline{D}_{\dot{\beta}_{T^*}} \dot{\beta}_{T^*} = 0, \tag{5}$$

where $\overline{D}$ is the connection on the hyperbolic unit sphere $S_1^2$ and the equation of tangent indicatrix is $\beta_{T^*} = T^*$. Thus from the Gauss equation we can write

$$D_{\dot{\beta}_{T^*}} \dot{\beta}_{T^*} = \overline{D}_{\dot{\beta}_{T^*}} \dot{\beta}_{T^*} + \varepsilon g\left(S(\dot{\beta}_{T^*}), \dot{\beta}_{T^*}\right) T^*, \varepsilon = g(T^*, T^*) = 1.$$

On the other hand, from the Lemma 2.2. straightforward computation gives

$$\dot{\beta}_{T^*} = t_{T^*} = \frac{d\beta_{T^*}}{ds} \frac{ds}{ds_{T^*}} = (\kappa \sinh\theta + \tau \cosh\theta) N \frac{ds}{ds_{T^*}}.$$

Moreover, we get

$$\frac{ds}{ds_{T^*}} = \frac{1}{\kappa \sinh\theta + \tau \cosh\theta}, \quad t_{T^*} = N,$$

$$D_{t_{T^*}} t_{T^*} = \frac{\kappa}{\kappa \sinh\theta + \tau \cos h\theta} T - \frac{\tau}{\kappa \sinh\theta + \tau \cos h\theta} B$$

and $g\left(S(t_{T^*}), t_{T^*}\right) = -1$.

Using these in the Gauss equation, we immediately have

$$\overline{D}_{t_{T^*}} t_{T^*} = \frac{\kappa}{\kappa \sinh\theta + \tau \cosh\theta} T - \frac{\tau}{\kappa \sinh\theta + \tau \cosh\theta} B + T^*.$$

From the Eq. (5) and Lemma 2.2. we get

$$\left(\frac{\kappa}{\kappa \sinh\theta + \tau \cosh\theta} + \sinh\theta\right) T + \left(\frac{\tau}{\kappa \sinh\theta + \tau \cosh\theta} + \cosh\theta\right) B = 0.$$

Since $T, N, B$ are linearly independent, we have

$$\begin{cases} \dfrac{\kappa}{\kappa sinh\theta + \tau cosh\theta} + ch\theta = 0 \\ \dfrac{\tau}{\kappa sinh\theta + \tau cosh\theta} + sh\theta = 0 \end{cases}$$

it follows that,

$$\kappa cosh\theta + \tau sinh\theta = 0, \tag{6}$$

$$\frac{\tau}{\kappa} = -coth\theta. \tag{7}$$

So from the Eq. (7) and Remark 2.1. we can give the following result.

**Result 3.1.1:** Let $(\alpha,\beta)$ be a timelike-spacelike Bertrand couple. If $\alpha$ is a general helix, then the tangent indicatrix $\beta_{T^*}$ of $\beta$ is a geodesic on $S_1^2$.

Moreover from Lemma 2.1. and Result 3.1.1. we can give the following theorem to characterize the natural lift of the tangent indicatrix of $\beta$ without proof.

**Theorem 3.1.1:** Let $(\alpha,\beta)$ be a timelike-spacelike Bertrand couple. If $\alpha$ is a general helix, then the natural lift $\overline{\beta}_{T^*}$ of the tangent indicatrix $\beta_{T^*}$ of $\beta$ is an integral curve of the geodesic spray on the tangent bundle $T(S_1^2)$.

## 3.2. The natural lift of the spherical indicatrix of the principal normal vectors of β

Let $\beta_{N^*}$ be the spherical indicatrix of principal normal vectors of $\beta$ and $\overline{\beta}_{N^*}$ be the natural lift of the curve $\beta_{N^*}$. If $\overline{\beta}_{N^*}$ is an integral curve of the geodesic spray, then by means of Lemma 2.1. we get,

$$\overline{D}_{t_{N^*}} t_{N^*} = 0 \tag{8}$$

that is

$$D_{t_{N^*}} t_{N^*} = \overline{D}_{t_{N^*}} t_{N^*} + \varepsilon g\left(S(t_{N^*}), t_{N^*}\right) N^* \quad \varepsilon = g(N^*, N^*) = 1$$

On the other hand, from Lemma 2.2. and Case I. straightforward computation gives

$$\dot{\beta}_{N^*} = t_{N^*} = cosh\varphi T - sinh\varphi B.$$

Moreover we get

$$D_{t_{N^*}} t_{N^*} = \frac{\varphi' cosh\varphi}{\|W\|} T + \frac{\kappa sinh\varphi - \tau cosh\varphi}{\|W\|} N - \frac{\varphi' sinh\varphi}{\|W\|} B \text{ and } g\left(S(t_{N^*}), t_{N^*}\right) = -1.$$

Using these in the Gauss equation, we immediately have

$$\overline{D}_{t_{N^*}} t_{N^*} = \frac{\varphi' \cosh\varphi}{\|W\|} T - \frac{\varphi' \sinh\varphi}{\|W\|} B.$$

Since *T, N, B* are linearly independent, we have

$$\begin{cases} \dfrac{\varphi' \cosh\varphi}{\|W\|} = 0 \\ \dfrac{\varphi' \sinh\varphi}{\|W\|} = 0 \end{cases},$$

it follows that,

$$\varphi' = 0, \tag{9}$$

$$\frac{\tau}{\kappa} = constant. \tag{10}$$

With similar computations for Case II. we have analogous result.

So from the Eq. (10) and Remark 2.1. we can give the following result.

**Result 3.2.1:** Let $(\alpha,\beta)$ be a timelike-spacelike Bertrand couple. If $\alpha$ is a general helix, then the principal normal indicatrix $\beta_{N^*}$ of $\beta$ is a geodesic on $S_1^2$.

Moreover from Lemma 2.1. and Result 3.2.1. we can give the following theorem to characterize the natural lift of the principal normal indicatrix of $\beta$ without proof.

**Theorem 3.2.1:** Let $(\alpha,\beta)$ be a timelike-spacelike Bertrand couple. If $\alpha$ is a general helix, then the natural lift $\overline{\beta}_{N^*}$ of the principal normal indicatrix $\beta_{N^*}$ of $\beta$ is an integral curve of the geodesic spray on the tangent bundle $T(S_1^2)$.

### 3.3. The natural lift of the spherical indicatrix of the binormal vectors of β

Let $\beta_{B^*}$ be the spherical indicatrix of binormal vectors of $\beta$ and $\overline{\beta}_{B^*}$ be the natural lift of the curve $\beta_{B^*}$. If $\overline{\beta}_{B^*}$ is an integral curve of the geodesic spray, then by means of Lemma 2.1. we get,

$$\overline{\overline{D}}_{t_{B^*}} t_{B^*} = 0 \tag{11}$$

that is

$$D_{t_{B^*}} t_{B^*} = \overline{\overline{D}}_{t_{B^*}} t_{B^*} + \varepsilon g\left(S(t_{B^*}), t_{B^*}\right) B^*, \quad \varepsilon = g(B^*, B^*) = -1,$$

On the other hand, from Lemma 2.2. straightforward computation gives

$$t_{B^*} = (\kappa\cosh\theta + \tau\sinh\theta) N \frac{ds}{ds_{B^*}}.$$

Moreover we get

$$\frac{ds}{ds_{B^*}} = \frac{1}{\kappa\cosh\theta + \tau\sinh\theta}, \quad t_{B^*} = N,$$

$$D_{t_{B^*}} t_{B^*} = \frac{\kappa}{\kappa\cosh\theta + \tau\sinh\theta} T - \frac{\tau}{\kappa\cosh\theta + \tau\sinh\theta} B \text{ and } g\left(S(t_{B^*}), t_{B^*}\right) = -1.$$

Using these in the Gauss equation, we immediately have

$$\overline{D}_{t_{B^*}} t_{B^*} = \frac{\kappa}{\kappa\cosh\theta + \tau\sinh\theta} T - \frac{\tau}{\kappa\cosh\theta + \tau\sinh\theta} B - B^*$$

From the Eq. (11) and Lemma 2.2. we get

$$\left(\frac{\kappa}{\kappa\cosh\theta + \tau\sinh\theta} - \cosh\theta\right) T + \left(-\frac{\tau}{\kappa\cosh\theta + \tau\sinh\theta} + \sinh\theta\right) B = 0.$$

Since $T, N, B$ are linearly independent, we have

$$\begin{cases} \dfrac{\kappa}{\kappa\cosh\theta + \tau\sinh\theta} - \cosh\theta = 0 \\ -\dfrac{\tau}{\kappa\cosh\theta + \tau\sinh\theta} + \sinh\theta = 0 \end{cases}$$

it follows that,

$$\kappa\sinh\theta + \tau\cosh\theta = 0, \qquad (12)$$

$$\frac{\tau}{\kappa} = \tanh\theta. \qquad (13)$$

So from the Eq. (13) and Remark 2.1. we can give the following result.

**Result 3.3.1:** Let $(\alpha, \beta)$ be a timelike-spacelike Bertrand couple. If $\alpha$ is a general helix, then the binormal indicatrix $\beta_{B^*}$ of $\beta$ is a geodesic on $H_0^2$.

Moreover from Lemma 2.1. and Result 3.3.1. we can give the following theorem to characterize the natural lift of the binormal indicatrix of $\beta$ without proof.

**Theorem 3.3.1:** Let $(\alpha, \beta)$ be a timelike Bertrand couple. If $\alpha$ is a general helix, then the natural lift $\overline{\beta}_{B^*}$ of the principal normal indicatrix $\beta_{B^*}$ of $\beta$ is an integral curve of the geodesic spray on the tangent bundle $T(H_0^2)$.

From the classification of all W-curves (i.e. a curves for which a curvature and a torsion are constants) in (Walrawe, 1995), we have following result with relation to curve $\alpha$.

**Result 3.3.2:**

(1) If the curve $\alpha$ with $\tau = 0$ then $\alpha$ is a planar timelike curve,

(2) If the curve $\alpha$ with $\tau = 0$, and $\kappa = constant > 0$ then $\alpha$ is a part of a orthogonal hyperbola,

(3) If the curve $\alpha$ with $\kappa = constant > 0$, $\tau = constant \neq 0$ and $|\tau| > \kappa$ then $\alpha$ is a part of a timelike circular helix,

$$\alpha(s) = \frac{1}{K}\left(\sqrt{\tau^2 K}\, s, \kappa\cos\left(\sqrt{K}\, s\right), \kappa\sin\left(\sqrt{K}\, s\right)\right) \text{ with } K = \tau^2 - \kappa^2.$$

(4) If the curve $\alpha$ with $\kappa = constant > 0$, $\tau = constant \neq 0$ and $|\tau| < \kappa$ then $\alpha$ is a part of a timelike hyperbolic helix,

$$\alpha(s) = \frac{1}{K}\left(\kappa\sinh\left(\sqrt{K}\, s\right), \left(\sqrt{\tau^2 K}\, s\right), \kappa\cosh\left(\sqrt{K}\, s\right)\right) \text{ with } K = \kappa^2 - \tau^2$$

From Lemma 3.1 in Choi et al 2012, we can write the following result:

**Result 3.3.3:** There is no timelike general helix with condition $|\tau| = |\kappa|$.

**Example:**

Let $\alpha(s) = \left(\frac{2\sqrt{3}}{3}s, \frac{1}{3}\cos\left(\sqrt{3}\, s\right), \frac{1}{3}\sin\left(\sqrt{3}\, s\right)\right)$ be a unit speed timelike circular helix with

$$\begin{cases} T = \left(\frac{2\sqrt{3}}{3}, \frac{\sqrt{3}}{3}\sin\left(\sqrt{3}\, s\right), \frac{\sqrt{3}}{3}\cos\left(\sqrt{3}\, s\right)\right) \\ N = \left(0, -\cos\left(\sqrt{3}\, s\right), -\sin\left(\sqrt{3}\, s\right)\right) \\ B = \left(-\frac{\sqrt{3}}{3}, \frac{2\sqrt{3}}{3}\sin\left(\sqrt{3}\, s\right), -\frac{2\sqrt{3}}{3}\cos\left(\sqrt{3}\, s\right)\right) \end{cases}, \kappa = 1 \text{ and } \tau = 2.$$

In this situation, spacelike Bertrand mate $\beta$ for $\alpha$ can be given by the equation

$$\beta(s) = \left(\frac{2\sqrt{3}}{3}s, \left(\frac{1}{3} - \lambda\right)\cos\left(\sqrt{3}\, s\right), \left(\frac{1}{3} - \lambda\right)\sin\left(\sqrt{3}\, s\right)\right), \lambda \in \mathbb{R}.$$

For $\lambda = \frac{4}{3}$, we have

$$\beta(s) = \left(\frac{2\sqrt{3}}{3}s, -\cos\left(\sqrt{3}\, s\right), -\sin\left(\sqrt{3}\, s\right)\right).$$

The straight forward calculations give the following spherical indicatrices and natural lift curves of spherical indicatrices for $\beta$,

$$\begin{cases} \beta_{T^*} = \left(\dfrac{2}{\sqrt{5}}, \dfrac{3}{\sqrt{5}}\sin(\sqrt{3}\,s), -\dfrac{3}{\sqrt{5}}\cos(\sqrt{3}\,s)\right) \\ \beta_{N^*} = \left(0, \cos(\sqrt{3}\,s), \sin(\sqrt{3}\,s)\right) \\ \beta_{B^*} = \left(\dfrac{3}{\sqrt{5}}, \dfrac{2}{\sqrt{5}}\sin(\sqrt{3}\,s), -\dfrac{2}{\sqrt{5}}\cos(\sqrt{3}\,s)\right) \end{cases},$$

$$\begin{cases} \overline{\beta}_{T^*} = \left(0, \dfrac{3\sqrt{3}}{\sqrt{5}}\cos(\sqrt{3}\,s), \dfrac{3\sqrt{3}}{\sqrt{5}}\sin(\sqrt{3}\,s)\right) \\ \overline{\beta}_{N^*} = \left(0, -\sqrt{3}\sin(\sqrt{3}\,s), \sqrt{3}\cos(\sqrt{3}\,s)\right) \\ \overline{\beta}_{B^*} = \left(0, \dfrac{2\sqrt{3}}{\sqrt{5}}\cos(\sqrt{3}\,s), \dfrac{2\sqrt{3}}{\sqrt{5}}\sin(\sqrt{3}\,s)\right) \end{cases},$$

respectively, (Figs. 2-7).

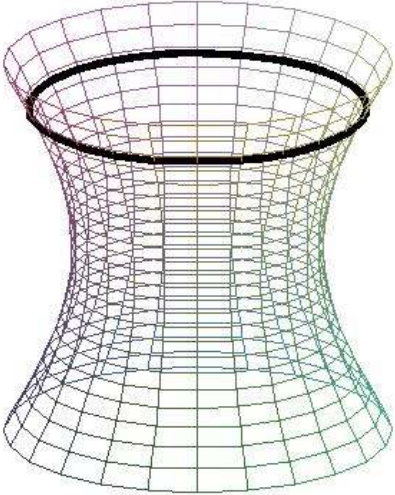
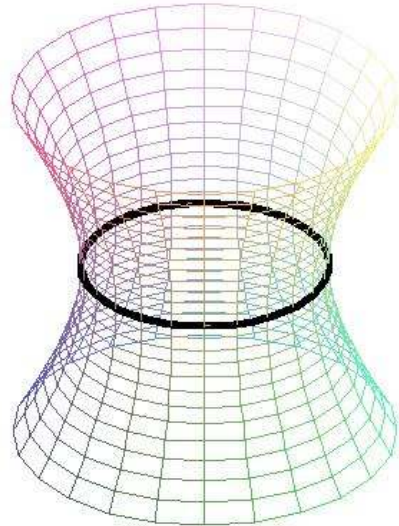

**Figure 2.** Tangent indicatrix for Bertrand mate of $\alpha$ on $S_1^2$

**Figure 3.** Principal normal indicatrix for Bertrand mate of $\alpha$ on $S_1^2$

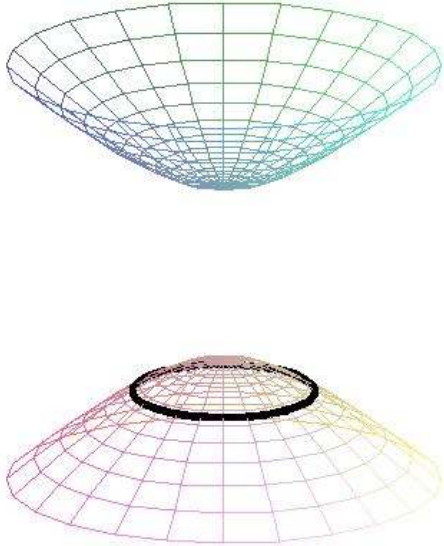
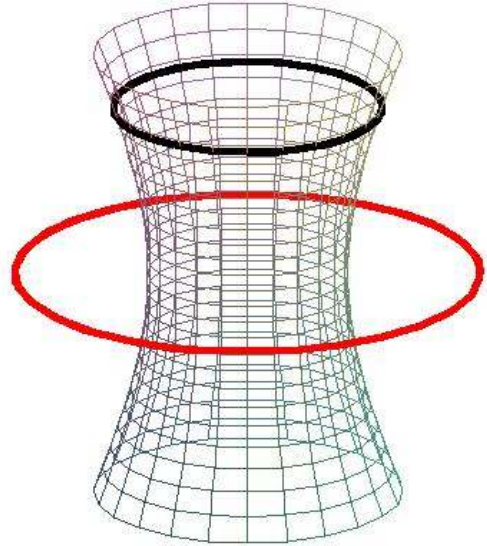

**Figure 4.** Binormal indicatrix for Bertrand mate of $\alpha$ on $H_0^2$

**Figure 5.** Tangent indicatrix for Bertrand mate of $\alpha$ and its natural lift

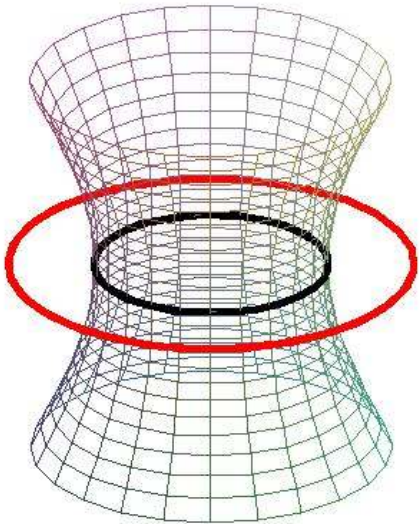
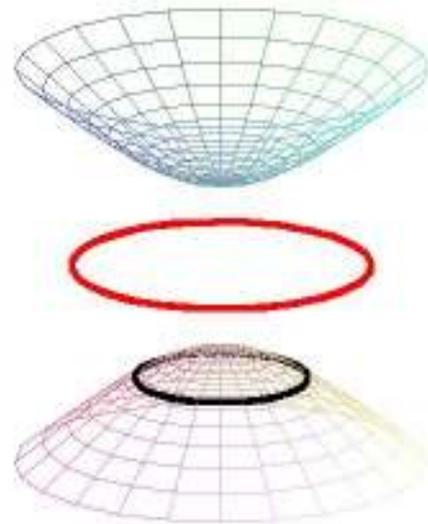

**Figure 6.** Principal normal indicatrix for Bertrand mate of $\alpha$ and its natural lift

**Figure 7.** Binormal indicatrix for Bertrand mate of $\alpha$ and its natural lift

**Discussion.** In this research, we extended the spherical indicatrix concept to the spacelike Bertrand mate of a timelike curve in the Minkowski 3-space $\mathbb{R}_1^3$. We investigated the natural lift curves and the geodesic sprays of the spherical indicatrices of tangent, principal normal and binormal vectors of the spacelike Bertrand mate $\beta$ of a timelike curve $\alpha$ and observed that $\alpha$ must be a general helix. It's expected that these results will be helpful to mathematicians who are specialized on mathematical modeling.

**Conflict of Interests**

The authors declare that there is no conflict of interests regarding the publication of this article.